\documentclass{article}
\usepackage{amssymb,amsfonts,amsmath,amsthm}
\usepackage{cleveref}
\usepackage{hyperref}
\usepackage{doi}
\usepackage{tikz}
\usepackage{sectsty}
\sectionfont{\centering}
\usepackage{OldStandard}
\usepackage[LGRgreek]{mathastext}
\usepackage[T1]{fontenc}
\usepackage{setspace}
\usepackage{fancyhdr}
\usepackage[bottom]{footmisc}
\usepackage[utf8]{inputenc}
\usepackage[english]{babel}
\usepackage{stmaryrd}
\usepackage{geometry}
\usepackage{blindtext}
\usepackage{fancyhdr}
\usepackage[all]{xy}

\newtheorem{theorem}{Theorem}
\newtheorem{definition}{Definition}
\newtheorem{proposition}{Proposition}
\usepackage{authblk}
\title{\textsc{Galois connections and isomorphism of simultaneous
ordered relations}}
\author{\textit{By}\textsc{ O. V. Atale }\footnote{Khandesh College Education Society's
Moolji Jaitha College Jalgaon-425001, Maharashtra, India. E-mail: atale.om@outlook.com}}
\affil{\textit{School of Mathematical Sciences, Moolji Jaitha College}\\ \textit{Jalgaon-425001, India}}
{\date{\textit{[Ann. of Comm. in Math., Vol. 6 (2) (2023), 109-117]}}}

\begin{document}
\maketitle

\begin{abstract}\fontsize{11}{12}\selectfont
In order theory, partially ordered sets are only equipped with one relation which decides the entire structure/Hasse diagram of the set. In this paper, we have presented how partially ordered sets can be studied under simultaneous partially ordered relations which we have called binary posets. The paper is motivated by the problem of operating a set simultaneously under two distinct partially ordered relations. It has been shown that binary posets follow the duality principle just like posets do. Within this framework, some new definitions concerning maximal and minimal elements are also presented. Furthermore, some theorems on order isomorphism and Galois connections are derived.
\end{abstract}

\fontsize{11}{12}\selectfont
\textbf{\textsc{\S 1. Introduction.}} Let $P$ be a non-empty set. A relation $\,\preceq$ on $P$ is known as partially ordered relation (POR) if $\,\preceq$ is \cite{1} :-
\begin{enumerate}
  \item Reflexive: $a\,\preceq{a}, \forall\, a\in P$.
  \item Anti-symmetric: $[(a\,\preceq\, b) \land (b\,\preceq a) ]\implies a=b, \forall\, a,b\in P$. 
  \item Transitive: $[(a\,\preceq\, b) \land (b\,\preceq \,c)] \implies a\,\preceq\, c, \forall\, a,b,c \in P$.
\end{enumerate}

In order theory, if $\,\preceq$ is a partially ordered relation on a non-empty set $P$, then an ordered pair $(P,\,\preceq)$ is known as a partially ordered set or poset. In this paper, we are going to construct the theory of posets equipped with binary relation ${\diamondsuit} = \left( {{ \,\preceq _1\,},{ \,\preceq _2\,}} \right)$. The partially ordered relations $\,\preceq_{1}$ and $\,\preceq_{2}$ may or may not be same. The reason for constructing the concept of a poset with binary relations is that it would help us to study a partially ordered set with two partially ordered relations simultaneously. Generally, when we study partially ordered sets, we take into account 3 properties: reflexivity, anti-symmetric and transitive property. Under the definition of these properties, we show under what conditions involving a single relation $\,\preceq$ can two different elements can be related to itself, equal to each other, and related to other elements. The drawback, however, is that this only helps us study non-sets under one relation at a time.\vspace{5mm}

\textbf{\textsc{\S 2. Binary poset.}} In this paper, we are going to modify the definition of properties of partially ordered relations so that we could study non-empty sets under two relations simultaneously at a time. We call this newly constructed poset a binary poset. It is shown that binary posts satisfy all the properties that a poset does. We can define isomorphism from one binary poset to another binary poset using the definition of a isotone, the binary posets also follow the duality principle. The derived theorems and definitions are then applied to derive some results on the power set. In the last section, we have given how a Galois connection between two binary posets can be defined. Some related examples and theorems on Galois's connections of binary posets are derived. Throughout the paper $\land$ and $\lor$ are short-hand mathematical notations for "and" and "or" respectively. 
 
\begin{definition}[partially ordered binary relation]
Let $P$ be a non-empty set. A relation ${\diamondsuit} = \left( {{ \,\preceq _1\,},{ \,\preceq _2\,}} \right)$ is known as a partially ordered binary relation (POBR) if ${\diamondsuit}$ is:-
\begin{enumerate}
  \item Reflexive: $a { \,\preceq _1\,}a  \,\preceq_2\, a ,\forall\, a  \in P.$
  \item Anti-symmetric: $\forall\,\, a ,b ,c  \in P$ $$\left[(a { \,\preceq _1\,}b { \,\preceq _2\,}c) \land (b { \,\preceq _1\,}a { \,\preceq _2\,}c) \land (a { \,\preceq _1\,}c { \,\preceq _2\,}b)\right] \Rightarrow (a  = b  = c).$$
  \item \textit{Transitive: $\forall\,\, a ,b, c ,d ,e  \in P$ $$\left[(a\, { \,\preceq _1\,}b { \,\preceq _2\,}c) \land (b { \,\preceq _1\,}d{\,\preceq_{2}}c) \land (c { \,\preceq _2\,}e) \right] \implies \left[ (a{\,\preceq _1\,}d { \,\preceq _2\,}c) \land (a\,\preceq_{1}b\,\preceq_{2}e)\right].$$}
\end{enumerate}
\end{definition}
Above properties are the extension of the binary relations from one that involves one operation to two operations. In the similar manner as above, one can define binary relations for other properties that include symmetric, connected, well-founded, 	irreflexive, and asymmetric properties of sets. Since, in this paper, we are only going to deal with partially ordered sets, we won't be needing them here. But in the case of reflexive, anti-symmetric, and transitive properties, the definition provided in Def. 1 can be thought of as universal relations operating on two partially ordered relations simultaneously.

\begin{definition}[Binary partially ordered set]
Let ${\diamondsuit} = \left( {{ \,\preceq _1\,},{ \,\preceq _2\,}} \right)$ be a partially ordered binary relation on a non-empty set $P$, then an ordered pair $\left( {P,\diamondsuit } \right)$ is known as binary partially ordered set (binary poset).
\end{definition}

The set of all natural numbers under relation $\diamondsuit=(\leq,|)$ forms a binary poset.

\begin{theorem}
If $\left( {P,\diamondsuit '} \right)$ and $\left( {P,\tilde \diamondsuit } \right)$ are two binary posets. Then $\left( {P,\diamondsuit ' \cap \tilde \diamondsuit } \right)$ is also a binary poset, where ${\diamondsuit ' = \left( {{ \,\preceq _1\,},{ \,\preceq _2\,}} \right)}$ and ${\tilde \diamondsuit  = \left( {{ \,\preceq _3\,},{ \,\preceq _4\,}} \right)}$.
\end{theorem}

\textit{Proof:} Suppose that $\left( {P,\diamondsuit '} \right)$ and $\left( {P,\tilde \diamondsuit } \right)$ are two binary posets. Let 
\begin{equation*}
a\,{{\preceq '}}\,b \,{{\preceq ''}}\,c  \iff \left[\left(a { \,\preceq _1\,}b { \,\preceq _2\,}c \right) \land \left(a { \,\preceq _3\,}b { \,\preceq _4\,}c \right)\right] \tag{2.1}  
\end{equation*} 
i.e. $\diamondsuit  = \diamondsuit ' \cap \tilde \diamondsuit $ where $a ,b ,c  \in P$. First, we prove reflexive property. Let $a\in P$. Since ${\diamondsuit '}$ and $\tilde \diamondsuit $ are reflexive, we get
\begin{equation*}
    \left[\left(a { \,\preceq _1\,}a { \,\preceq _2\,}a\right) \land \left(a { \,\preceq _3\,}a { \,\preceq _4\,}a\right)\right] \implies (a {{ \,\preceq '}}\,a {{ \,\preceq ''}}\,a).  \tag{2.2}
\end{equation*} 
Hence, $\diamondsuit $ is reflexive. Now, we prove the anti-symmetric property. Let $a,b,c\in P$. Since ${\diamondsuit '}$ and $\tilde \diamondsuit $ are anti symmetric, we get 
\begin{equation*}
 \left[(a { \,\preceq _1\,}b { \,\preceq _2\,}c) \land (b { \,\preceq _1\,}a { \,\preceq _2\,}c) \land (a { \,\preceq _1\,}c { \,\preceq _2\,}b)\right],\tag{2.3}   
\end{equation*}
\begin{equation*}
\left[(a { \,\preceq _3\,}b { \,\preceq _4\,}c) \land (b { \,\preceq _3\,}a { \,\preceq _4\,}c) \land (a { \,\preceq _3\,}c { \,\preceq _4\,}b)\right].\tag{2.4}   
\end{equation*}
Thus 
\begin{equation*}
 \left[(a { { \,\preceq '}}b { { \,\preceq ''}}c) \land (b { { \,\preceq '}}a { { \,\preceq ''}}c) \land (a { { \,\preceq '}}c { { \,\preceq ''}}b)\right]\implies (a=b=c).\tag{2.5}   
\end{equation*}
Hence, $\diamondsuit$ is anti-symmetric. Now, we prove transitive property. Let $a,b,c,d,e  \in P$. Since ${\diamondsuit '}$ and $\tilde \diamondsuit $ are transitive, we get
\begin{equation*}
 \left[(a { \,\preceq _1\,}b { \,\preceq _2\,}c) \land (b { \,\preceq _1\,}d{ \,\preceq _2\,}c) \land (c{ \,\preceq _2\,}e)\right]  \implies \left[ (a { \,\preceq _1\,}d { \,\preceq _2\,}c) \land (a { \,\preceq _1\,}b { \,\preceq _2\,}e)\right],\tag{2.6}   
\end{equation*}
\begin{equation*}
  \left[(a { \,\preceq _3\,}b { \,\preceq _4\,}c) \land (b { \,\preceq _3\,}d{ \,\preceq _4\,}c) \land (c{ \,\preceq _4\,}e)\right]  \implies \left[ (a { \,\preceq _3\,}d { \,\preceq _4\,}c) \land (a { \,\preceq _3\,}b { \,\preceq _4\,}e)\right].\tag{2.7}  
\end{equation*}
Thus, 
\begin{equation*}
 \left[ (a { { \,\preceq '}}b { { \,\preceq ''}}c) \land (b { { \,\preceq '}}d{\,\preceq ''}c) \land (c{ { \,\preceq ''}}e)\right]  \implies \left[(a { { \,\preceq '}}d{ { \,\preceq ''}}c) \land (a { { \,\preceq '}}b { { \,\preceq ''}}e)\right].\tag{2.8}   
\end{equation*}
Hence, $\diamondsuit$ is transitive. Therefore, $(P,\diamondsuit)=\left( {P,\diamondsuit ' \cap \tilde \diamondsuit } \right)$ is a binary poset. The above theorem can be generalized as follows.

\begin{theorem}
If $\left( {P,{\diamondsuit _n}} \right), n = 0,1,2,...,m$ are binary posets where ${\diamondsuit _n} = \left( {{ \,\preceq _{n + 1}},{ \,\preceq _{n + 2}}} \right)$. Then, $\left( {P,\diamondsuit } \right)$ is also a binary poset where $\diamondsuit  = {\diamondsuit _0} \cap {\diamondsuit _1\,} \cap {\diamondsuit _2\,} \cap ... \cap {\diamondsuit _m}$.
\end{theorem}\vspace{5mm}

\textbf{\textsc{\S3. Maximal and minimal elements.}} Let $(P,\preceq)$ be a poset. Then, the notation $a(\forall,\in,\notin)P$ means "\textit{for all, or for some, or no}"$a\in P$. For example, let $a,b\in P$ and $a\,\preceq\, b$. Now, $a$ being fixed, there might be a possibility that $a\,\preceq\, b$ will be true for all $b\in P$, or for some $b\in P$ or for no $b\in P$. Therefore we collectively write $a\,\preceq\, b,\,b\,(\forall,\in,\notin)\,P$
\begin{definition}[Minimal and maximal greatest element]
Let $(P,\diamondsuit)$ be a binary poset where $\diamondsuit=(\,\preceq_{1},\,\preceq_{2})$. Let $x, y\in P$ be any two elements such that 
\begin{equation*}
 \left[(a { \,\preceq _1\,}x{ \,\preceq _2\,}c),(\forall\, a\in P) \land c\, (\forall,\in,\notin)\,P\right] \land \left[(a { \,\preceq _1\,}b { \,\preceq _2\,}y), (\forall\, b\in P) \land a\,(\forall,\in,\notin)\,P\right]. \tag{3.1}  
\end{equation*}
Then, $g_{max} = \sup \left\{ {x,y} \right\}$ is known as the maximal greatest element and $g_{min}= \inf \left\{ {x,y} \right\}$ is known as the minimal greatest element.
\end{definition}

The above definition can be summarized as follows. Lets say we have a non-empty set $P$ equipped with partially ordered binary relation $\diamondsuit  = \left( {{ \,\preceq _1\,},{ \,\preceq _2\,}} \right)$. Then, the greatest elements generated by $\,\preceq_{1}$ and $\,\preceq_{2}$ may or may not be equal. If they are equal, then we are through. If they are not equal, then the greatest one will be the maximal greatest element, and the smallest one will be the minimal greatest element. Similarly, for the minimal and maximal least element, we have the following definition.

\begin{definition}[Minimal and maximal least element]
Let $(P,\diamondsuit)$ be a binary poset where $\diamondsuit=(\,\preceq_{1},\,\preceq_{2})$. let $x, y\in P$ be any two elements such that
\begin{equation*}
\left[(x { \,\preceq _1\,}b{ \,\preceq _2\,}c),(\forall\, b\in P) \land c\, (\forall,\in,\notin)\,P\right] \land \left[(a { \,\preceq _1\,}y { \,\preceq _2\,}c), (\forall\, c\in P) \land a\,(\forall,\in,\notin)\,P\right].    \tag{3.2}
\end{equation*}
Then, $l_{max} = \sup \left\{ {x,y} \right\}$ is known as the maximal least element and $l_{min}= \inf \left\{ {x,y} \right\}$ is known as the minimal least element.
\end{definition}

\begin{definition}[Bounded binary poset]
A binary poset is known as bounded if it has the maximal greatest element and minimal least element.
\end{definition}
\begin{definition}[Unbounded binary poset]
A binary poset is known as unbounded if it is not a bounded binary poset.
\end{definition}

\begin{theorem}
Let $(P,\diamondsuit)$ be a binary poset where $\diamondsuit=(\,\preceq_{1},\,\preceq_{2})$. Then, the maximal greatest element (if it exists) is unique.
\end{theorem}

\textit{Proof:} Suppose that $g'$ and $\tilde g$ are the two maximal greatest elements of a binary poset $(P,\diamondsuit)$ where $\diamondsuit=(\,\preceq_{1},\,\preceq_{2})$. Using definition 2, we can write $g' = \sup \left\{ {{g_1\,},{g_2\,}} \right\}$ and $\tilde g = \sup \left\{ {{g_3\,},{g_4\,}} \right\}$ where $g_{1},g_{2},g_{3},g_{4}\in P$ such that 
\begin{equation*}
\left[(a { \,\preceq _1\,}g_{1}{ \,\preceq _2\,}c),(\forall\, a\in P) \land c\, (\forall,\in,\notin)\,P\right] \land \left[(a { \,\preceq _1\,}b { \,\preceq _2\,}g_{2}), (\forall\, b\in P) \land a\,(\forall,\in,\notin)\,P\right].    \tag{3.3}
\end{equation*}
\begin{equation*}
 \left[(g_{3} { \,\preceq _1\,}g_{1}{ \,\preceq _2\,}c), c\,(\forall,\in,\notin)\,P\right] \land \left[(a { \,\preceq _1\,}g_{4} { \,\preceq _2\,}g_{2}), a\,(\forall,\in,\notin)\,P\right].   \tag{3.4}
\end{equation*}
Now, $g_{3}$, $g_{4}\in P$. Therefore,
\begin{equation*}
 \left[(g_{1} { \,\preceq _1\,}g_{3}{ \,\preceq _2\,}c), c\,(\forall,\in,\notin)\,P\right] \land \left[(a { \,\preceq _1\,}g_{2} { \,\preceq _2\,}g_{4}), a\,(\forall,\in,\notin)\,P\right].   \tag{3.5}
\end{equation*}
Also. $g_{1}$, $g_{2}\in P$. Therefore, Now, $g_{3}$, $g_{4}\in P$. This implies $g_{1}=g_{3}$ and $g_{2}=g_{4}$. And hence, $g' = \tilde g$. This completes our proof.

Similarly, for the uniqueness of the minimal greatest element we have the following theorem.

\begin{theorem}
Let $(P,\diamondsuit)$ be a binary poset where $\diamondsuit=(\,\preceq_{1},\,\preceq_{2})$. Then, the minimal greatest element (if it exists) is unique.
\end{theorem}
The same result applies to the uniqueness of minimal and maximal least elements.
\begin{theorem}
Let $(P,\diamondsuit)$ be a binary poset where $\diamondsuit=(\,\preceq_{1},\,\preceq_{2})$. Then, the maximal least element (if it exists) is unique.
\end{theorem}

\begin{theorem}
Let $(P,\diamondsuit)$ be a binary poset where $\diamondsuit=(\,\preceq_{1},\,\preceq_{2})$. Then, the minimal least element (if it exists) is unique.
\end{theorem}

Proof of Theorem 5 and 6 follows from the same technique of proof of Theorem 3. 

\begin{theorem}
Let $X$ be a non-empty set and $P(X)$ be the power set of X. Define $\diamondsuit  = \left( {{ \,\preceq _1\,},{ \,\preceq _2\,}} \right) = \left( { \subseteq , \subseteq } \right)$ on $P(X)$ by 
\begin{equation*}
  A{ \,\preceq _1\,}B{ \,\preceq _2\,}C \iff A \subseteq B \subseteq C.  \tag{3.6}
\end{equation*}
Then $(P(X),\diamondsuit)$ is a binary poset.
\end{theorem}

\textit{Proof:} Let $A,B,C,D,E\in P(X)$. Since every subset is a subset of itself, we can show that $A\subseteq A\subseteq A, \forall\, A\in P(X)$. Therefore, $\diamondsuit$ is reflexive. Secondly, since 
\begin{align*}
 \big[(A&\subseteq B\subseteq C) \land (B\subseteq A\subseteq C)\land (A\subseteq C\subseteq B)\big] \\&\implies (A=B=C), \forall\, A,B,C\in P(X), \tag{3.7}   
\end{align*}
therefore $\diamondsuit$ is anti-symmetric. Now, we know that if $A\subseteq B$ and $B\subseteq C$, then $A\subseteq C$. Therefore, if
\begin{equation*}
  \left[(A\subseteq B\subseteq C) \land (B\subseteq D \subseteq C) \land (C\subseteq E)\right] \implies \left(A\subseteq D\subseteq C)\land (A\subseteq B\subseteq E)\right].\tag{3.8}  
\end{equation*}
Hence, $\diamondsuit$ is a partial order binary relation and thus $(P(X),\diamondsuit)$ is a binary poset.\vspace{5mm}

\textbf{\textsc{\S4. Binary poset isomorphism and duality principle.}} In this section, we will show that binary posets also follow the duality principle.

\begin{definition}[Binary poset isomorphism]
Let $(P,\diamondsuit)$ and $(Q,\tilde\diamondsuit)$ be two binary posets where $\diamondsuit  = \left( {{ \,\preceq _1\,},{ \,\preceq _2\,}} \right)$ and $\tilde\diamondsuit  = \left( {{ \,\preceq _3\,},{ \,\preceq _4\,}} \right)$, A function $\psi:(P,\diamondsuit)\to (Q,\tilde{\diamondsuit})$ is known as isomorphism if
\begin{enumerate}
    \item For $a,b,c\in P$, $a { \,\preceq _1\,}b { \,\preceq _2\,}c  \iff \psi \left( a  \right){ \,\preceq _3\,}\psi \left( b  \right){ \,\preceq _4\,}\psi \left( c  \right)$, 
    \item $\psi$ is bijective.
\end{enumerate}
Symbolically, we write $(P,\diamondsuit)\cong (Q,\tilde{\diamondsuit})$ 
\end{definition}
\begin{definition}[Isotone]
A function $\psi :\left( {P,\diamondsuit } \right) \to \left( {Q,\diamondsuit } \right)$ where $\diamondsuit  = \left( {{ \,\preceq _1\,},{ \,\preceq _2\,}} \right)$ is known as isotone if
\begin{equation*}
 a { \,\preceq _1\,}b { \,\preceq _2\,}c  \iff \psi \left( c  \right){ \,\preceq _1\,}\psi \left( b  \right){ \,\preceq _2\,}\psi \left( c  \right), \forall\, a,b,c\in P. \tag{4.1}  
\end{equation*}
\end{definition}
\begin{theorem}
A function $\psi :\left( {P,\diamondsuit } \right) \to \left( {Q,\tilde \diamondsuit } \right)$ is a binary poset isomorphism iff $\psi$ is isotone and has isotone inverse. 
\end{theorem}

\textit{Proof:} Suppose that $\psi :\left( {P,\diamondsuit } \right) \to \left( {Q,\diamondsuit } \right)$ is a binary poset isomorphism where $\diamondsuit  = \left( {{ \,\preceq _1\,},{ \,\preceq _2\,}} \right)$. Therefore, $\psi$ is bijective and satisfis the condition that for $a,b,c\in P$,
\begin{equation*}
(a { \,\preceq _1\,}b { \,\preceq _2\,}c)  \iff \psi \left( c  \right){ \,\preceq _1\,}\psi \left( b  \right){ \,\preceq _2\,}\psi \left( c  \right).\tag{4.2}   
\end{equation*}
By definition of isotone, $\psi$ is an isotone if it is bijective. Therefore, ${\psi ^{ - 1}}:\left( {P,\diamondsuit } \right) \to \left( {Q,\diamondsuit } \right)$ exists. Let $\delta,\epsilon,\zeta\in Q$ such that $\delta { \,\preceq _1\,}\epsilon { \,\preceq _2\,}\zeta $. Since $\psi$ is surjective, $\exists\, a ,b ,c  \in P$ such that $\psi \left( a  \right) = \delta $ and $\psi \left( b  \right) = \epsilon ,\psi \left( c  \right) = \zeta $. Therefore $a  = \psi^{-1}\left( \delta  \right),b  = {\psi ^{ - 1}}\left( \epsilon  \right)$, $c  = {\psi ^{ - 1}}\left( \zeta  \right)$. Now, 
\begin{align*}
(\delta { \,\preceq _1\,}\epsilon { \,\preceq _2\,}\zeta)  &\Rightarrow \psi \left( a  \right){ \,\preceq _1\,}\psi \left( b  \right){ \,\preceq _2\,}\psi \left( c  \right) \\&\Rightarrow (a { \,\preceq _1\,}b { \,\preceq _2\,}c)  \\&\Rightarrow {\psi ^{ - 1}}\left( \delta  \right){ \,\preceq _1\,}{\psi ^{ - 1}}\left( \epsilon  \right){ \,\preceq _2\,}{\psi ^{ - 1}}\left( \zeta  \right). \tag{4.3}  
\end{align*}
Therefore, $\psi^{-1}$ is isotone.

\begin{definition}[Dual]
Let $\diamondsuit = \left( {{ \,\preceq _1\,},{ \,\preceq _2\,}} \right)$ be a relation on the set $P$. Then, the dual of $\diamondsuit$ is denoted by $\tilde{\diamondsuit}= \left( {{{\tilde\preceq }_1\,},{{\tilde\preceq }_2\,}} \right)$ and is defined as $a {{\tilde\preceq }_1\,}b {{\tilde\preceq }_2\,}c  \iff c {  \,\preceq }_2\,b {{ \,\preceq }_1\,}a$ or $a {{\tilde\preceq }_1\,}b {{\tilde\preceq }_2\,}c  \iff b {{  \,\preceq }_1\,}a$ and $c {{  \,\preceq }_2\,}b, \forall\, a,b,c\in P$.
\end{definition}

The two different definitions of dual provided above are equal whereas the second definition is easy to use in the proof as compared to the first.

\begin{theorem}[Duality principle for binary posets]
Let $(P,\diamondsuit)$ be a binary poset, then $(P,\tilde{\diamondsuit})$ is also a binary poset.
\end{theorem}

\textit{Proof:}  Let $\diamondsuit  = \left( {{ \,\preceq _1\,},{ \,\preceq _2\,}} \right)$, $\tilde \diamondsuit  = \left( {{{\tilde\preceq }_1\,},{{\tilde\preceq }_2\,}} \right)$ and $a,b,c\in P$. We know that $a { \,\preceq _1\,}a$ and $a{ \,\preceq _2\,}a $, thus ${a {{\tilde\preceq }_1\,}a {{\tilde\preceq }_2\,}a }$. Hence $\tilde\diamondsuit$ is reflexive. Let
\begin{equation*}
(a { \tilde\preceq _1\,}b { \tilde\preceq _2\,}c) \land (b { \tilde\preceq _1\,}a { \tilde\preceq _2\,}c).    \tag{4.4}
\end{equation*}
Therefore, Eqn. (4.4) implies
\begin{equation*}
(b {{  \,\preceq }_1\,}a ) \land (c {{  \,\preceq }_2\,}b) \land (a {\,\preceq _1\,}b ) \land (c\,\preceq_{2}a)\implies(a=b).    \tag{4.5}
\end{equation*}
Similarly, we can further prove that $a=b=c$. Thus, $\tilde\diamondsuit$ is anti-symmetric. Now, let 
\begin{equation*}
 (a {{\tilde\preceq }_1\,}b {{\tilde\preceq }_2\,}c) \land (b {{\tilde\preceq }_1\,}d {{\tilde\preceq }_2\,}c) \land (c {{\tilde\preceq }_2\,}e).   \tag{4.6}
\end{equation*}
Therefore, using Def. 9, 
\begin{equation*}
 \left[(b { \,\preceq _1\,}a) \land (c { \,\preceq _2\,}b) \land (d { \,\preceq _1\,}b) \land (c { \,\preceq _2\,}d) \land (e { \,\preceq _2\,}c)\right],   \tag{4.7}
\end{equation*}
\begin{equation*}
 \implies \left[(d { \,\preceq _1\,}a) \land (c { \,\preceq _2\,}d)\right]  \implies (a {{\tilde\preceq }_1\,}d {{\tilde\preceq }_2\,}c)   \tag{4.8}
\end{equation*}
and
\begin{equation*}
   \left[(b { \,\preceq _1\,}a) \land (e { \,\preceq _2\,}b)\right]  \implies (a {{\tilde\preceq }_1\,}b {{\tilde\preceq }_2\,}e). \tag{4.9}
\end{equation*}
Hence, $\diamondsuit$ is transitive and $(P,\tilde\diamondsuit)$ is also a binary poset.

\begin{definition}
Let $(P,\diamondsuit)$ be a binary poset. Then, $(P,\tilde\diamondsuit)$ is known as the dual of $(P,\diamondsuit)$.
\end{definition}

\begin{definition}[Self dual of a binary poset]
If $(P,\diamondsuit)$ is a binary poset and $(P,\diamondsuit) \cong (P,\tilde\diamondsuit)$, then $(P,\diamondsuit)$ is known as self dual poset.
\end{definition}

\begin{theorem}
Let $X$ be a non-empty set. Then, $(P(X),\diamondsuit)$ is a self dual poset, where $\diamondsuit=(\subseteq,\subseteq)$.
\end{theorem}

\textit{Proof:} Define the mapping $\psi :P\left( X \right) \to P\left( {\tilde X} \right)$ by $\psi \left( A \right) = X - A,\forall\, A \in P\left( X \right)$. Let $A,B\in P(X)$ such that $A=B$ where $P\left( {\tilde X} \right)=(P(X),\tilde\diamondsuit)$. Therefore, 
\begin{equation*}
 X-A=X-B\implies \psi(A)=\psi(B).   \tag{4.10}
\end{equation*}
Hence $\psi$ is well-defined and injective. Now, let $A\in P(\tilde{X})$. Therefore 
\begin{equation*}
 X-A\in P(X)\implies \psi(X-A)=X-(X-A)=B.   \tag{4.11}
\end{equation*}
Hence $\psi$ is surjective. Let $A,B,C\in P(X)$ such that $A\subseteq B\subseteq C$. Therefore, 
\begin{align*}
 A\subseteq B\subseteq C&\implies X-B\subseteq X-A\subseteq C \\& \implies \psi(B)\subseteq \psi(A)\subseteq C \\& \implies C \supseteq \psi \left( A \right) \supseteq \psi \left( B \right).\tag{4.12}   
\end{align*}
Also, 
\begin{align*}
    A\subseteq B\subseteq C&\implies A\subseteq X-C\subseteq X-B \\&\implies A\subseteq\psi(C)\subseteq\psi(B)\\&\implies \psi \left( B \right) \supseteq \psi \left( C \right) \supseteq A. \tag{4.13}
\end{align*}
Therefore
\begin{equation*}
 C \supseteq \psi \left( A \right) \supseteq \psi \left( B \right)\,\,\mathrm{and}\,\, \psi \left( B \right) \supseteq \psi \left( C \right) \supseteq A\implies \psi \left( A \right) \supseteq \psi \left( B \right) \supseteq \psi \left( C \right).   \tag{4.14}
\end{equation*}
Hence, $\psi :P\left( X \right) \to P\left( {\tilde X} \right)$ is an ispmorphism and $P\left( X \right) \cong P\left( {\tilde X} \right)$. Therefore, $(P(X)\diamondsuit)$ is self dual.

\begin{theorem}
Let $(P,\diamondsuit)$ be a binary poset and let $(P,\tilde{\tilde{\diamondsuit}})$ be the dual of $(P,\tilde{\diamondsuit})$. Then, $(P,\diamondsuit)$ is isomorphic to $(P,\tilde{\tilde{\diamondsuit}})$, i.e. $(P,\diamondsuit)\cong (P,\tilde{\tilde{\diamondsuit}})$. 
\end{theorem}

\textit{Proof:} Define $\psi:(P,\diamondsuit)\to(P,\tilde{\tilde{\diamondsuit}})$ by $\psi(a)=a, \forall\, a\in P$. Clearly, $\psi(a)$ is well defined and bijective. Let $\diamondsuit  = \left( {{ \,\preceq _1\,},{ \,\preceq _2\,}} \right),\tilde \diamondsuit  = \left( {{{\tilde\preceq }_1\,},{{\tilde\preceq }_2\,}} \right),\tilde{\tilde{\diamondsuit}}  = \left( {{{\tilde{\tilde{{\,\preceq }_1\,}}}},{{\tilde{\tilde{\,\preceq_2\,}}}}} \right)$ (by Def. 9). Now, Let $a,b,c \in X$ such that 
\begin{align*}
c { \,\preceq _2\,}b { \,\preceq _1\,}a  &\iff a\, {{\tilde\preceq }_1\,}b\, {{\tilde\preceq }_2\,}c \\&\iff c {{\tilde{\tilde{ \,\preceq_2\,}}}}b {{\tilde{\tilde{\,\preceq_1\,}}}}a \\& \iff \psi \left( c \right){{\tilde{\tilde{\preceq_2}}}}\,\psi \left( b  \right){{\tilde{\tilde{\preceq_1}}}}\,\psi \left( a \right).   \tag{4.15} 
\end{align*}
Therefore, $\psi(P,\diamondsuit)\to(P,\tilde{\tilde{\diamondsuit}})$ is an isomorphism and thus $(P,\diamondsuit)\cong (P,\tilde{\tilde{\diamondsuit}})$. \vspace{5mm}

\textbf{\textsc{\S5. Galois connections.}} Let $(P,\diamondsuit)$ be a binary poset where $\diamondsuit=(\,\preceq_{1},\,\preceq_{2})$ and $a,b\in P$. Then, define the notation $a\, \diamondsuit\, {b}=a\,\preceq_{1}b, a\,\preceq_{2}b$. Throughout the sequel we let $\diamondsuit=(\,\preceq_{1},\,\preceq_{2})$, $\tilde\diamondsuit=(\,\preceq_{3},\,\preceq_{4})$ and $\tilde{\tilde{\diamondsuit}}=(\,\preceq_{3},\,\preceq_{4})$.

\begin{definition}[Galois connection]
Let $(P,\diamondsuit)$ and $(Q,\tilde\diamondsuit)$ be two binary posets. Let $\psi^{*}:(P,\diamondsuit)\to (Q,\tilde\diamondsuit)$ and $\psi_{*}:(Q,\tilde\diamondsuit)\to (P,\diamondsuit)$ be a pair of functions. For $a\in P$ and $b\in Q$, if
\begin{equation*}
 \psi^{*}(a)\, \tilde\diamondsuit\, {b}\iff a\, \diamondsuit\, \psi_{*}(b),   \tag{5.1}
\end{equation*}
then the pair $\psi^{*}, \psi_{*}$ forms a Galois connection between binary posets $(P,\diamondsuit)$ and $(Q,\tilde\diamondsuit)$ and is denoted by $\mathrm{Gal\left(\psi^{*},\psi_{*}\right)}$ or more generally, $\mathrm{Gal\left(\psi^{*},\psi_{*}\right)}:(P,\diamondsuit)\to (Q,\tilde\diamondsuit)$.
\end{definition}

Following are some examples of a pair of binary posets and functions that form a Galois connection.

\begin{enumerate}
    \item  Let the map $\psi$ be an binary poset isomorphism between $(P,\diamondsuit)$ and $(Q,\tilde\diamondsuit)$. Then, $\mathrm{Gal\left(\psi,\psi^{-1}\right)}:(P,\diamondsuit)\to (Q,\tilde\diamondsuit)$ is a Galois connection.
    \item Let $(\mathbb{N},\diamondsuit)$ and $(\mathbb{Q}^{+},\diamondsuit)$ be to binary posets where $\diamondsuit=(\,\preceq,\,\preceq)$. Let $\psi^{*}:\mathbb{N}\to\mathbb{Q}^{+}$ be that standard embedding of the natural numbers into the rationals and $\psi_{*}:\mathbb{Q}^{+}\to\mathbb{N}$ be the map a positive rationals to the natural numbers corresponding to its integral part. Then, $\mathrm{Gal(\psi^{*},\psi_{*})}:(\mathbb{N},\diamondsuit)\to(\mathbb{Q}^{+},\diamondsuit)$ is a Galois connection.
    \item Let $(P,\diamondsuit)$ be an arbitrary binary poset and let $(\left\{0\right\},\tilde\diamondsuit)$ be a singleton binary poset where $\tilde\diamondsuit=(=,=)$. Let $\psi^{*}:(P,\diamondsuit)\to(\left\{0\right\},\tilde\diamondsuit)$ be a trivial function mapping all elements of $P$ to $\left\{0\right\}$ and $\psi_{*}:(\left\{0\right\},\tilde\diamondsuit)\to(P,\diamondsuit)$ be another dunction mapping $0$ to particular elements of $P$. Then, $\mathrm{Gal(\psi^{*},\psi_{*})}:(P,\diamondsuit)\to(\left\{0\right\},\tilde\diamondsuit)$ is a Galois connection.
\end{enumerate}

Following is one of the alternative definition of Def. 12.
\begin{theorem}
Let $(P,\diamondsuit)$ and $(Q,\tilde\diamondsuit)$ be two binary posets. Let $\psi^{*}:(P,\diamondsuit)\to (Q,\tilde\diamondsuit)$ and $\psi_{*}:(Q,\tilde\diamondsuit)\to (P,\diamondsuit)$ be a pair of functions. Then, $\mathrm{Gal\left(\psi^{*},\psi_{*}\right)}:(P,\diamondsuit)\to (Q,\tilde\diamondsuit)$ is Galois connection if:-
\begin{enumerate}
    \item $\psi^{*}$ and $\psi_{*}$ both are isotones, and
    \item $\forall\, a\in P$ and $b\in Q$, $a\, \diamondsuit\, \psi^{*}(\psi_{*}(a))$ and $\psi_{*}(\psi^{*}(b))\, \tilde\diamondsuit\, b $.
\end{enumerate}
\end{theorem}

\textit{Proof:} Suppose that $\mathrm{Gal\left(\psi^{*},\psi_{*}\right)}:(P,\diamondsuit)\to (Q,\tilde\diamondsuit)$ is Galois connection. Then, by definition we have 
\begin{equation*}
 \psi^{*}(a)\, \tilde{\diamondsuit}\, \psi^{*}(a) \iff a\, \diamondsuit\, \psi_{*}\circ\psi^{*}(a).   \tag{5.2}
\end{equation*}
Since $\tilde\diamondsuit$ is reflexive, $\psi^{*}(a)\, \tilde{\diamondsuit}\, \psi^{*}(a)$ holds true and thus $a\, \diamondsuit\, \psi_{*}\circ\psi^{*}(a)$. By a similar argument, we get $\psi_{*}\circ \psi^{*}(b)\, \tilde\diamondsuit\, b $. Now, let $a '\in P$ and $a\, \diamondsuit\, a'$. Since we have just shown that  $a'\, \diamondsuit\, \psi_{*}\circ\psi^{*}(a')$, we get $a\, \diamondsuit\, \psi_{*}\circ\psi^{*}(a')$. But, by Def. 12 we have 
\begin{equation*}
 \psi^{*}(a)\, \tilde\diamondsuit\, \psi^{*}(a ') \iff a\, \diamondsuit\, \psi_{*}\circ\psi^{*}(a '),   \tag{5.3}
\end{equation*}
therefore, $\psi^{*}\, \tilde\diamondsuit\, \psi^{*}(a,)$ and thus $\psi^{*}$ is a isotone. By a similar argument, we can show that $\psi_{*}$ is also a isotone.

So far, we are dealing with Galois connections that are equipped with binary partially ordered relations that are not same. When $\diamondsuit=\tilde\diamondsuit$, we can define monotone and antitone Galois connection as follows. 

\begin{definition}[Monotone Galois connection]
Let $(P,\diamondsuit)$ and $(Q,\diamondsuit)$ be two binary posets. Let $\psi^{*}:(P,\diamondsuit)\to (Q,\diamondsuit)$ and $\psi_{*}:(Q,\diamondsuit)\to (P,\diamondsuit)$ be a pair of functions. For $a\in P$ and $b\in Q$, if 
\begin{equation*}
 \psi^{*}(a)\, \diamondsuit\, {b} \iff a\, \diamondsuit\, \psi_{*}(b),   \tag{5.4}
\end{equation*}
then the pair $\psi^{*}, \psi_{*}$ forms a monotone Galois connection between binary posets $(P,\diamondsuit)$ and $(Q,\diamondsuit)$ and is denoted by $\mathrm{Gal_{mon}}\left(\psi^{*},\psi_{*}\right)$ or more generally, $\mathrm{Gal_{mon}}\left(\psi^{*},\psi_{*}\right):(P,\diamondsuit)\to (Q,\tilde\diamondsuit)$.
\end{definition}

\begin{definition}[Antitone Galois connection]
Let $(P,\diamondsuit)$ and $(Q,\diamondsuit)$ be two binary posets. Let $\psi^{*}:(P,\diamondsuit)\to (Q,\diamondsuit)$ and $\psi_{*}:(Q,\diamondsuit)\to (P,\diamondsuit)$ be a pair of functions. For $a\in P$ and $b\in Q$, if 
\begin{equation*}
  {b}\, \diamondsuit\, \psi^{*}(a) \iff a\, \diamondsuit\, \psi_{*}(b).  \tag{5.5}
\end{equation*}
Then the pair $\psi^{*}, \psi_{*}$ forms a antitone Galois connection between binary posets $(P,\diamondsuit)$ and $(Q,\diamondsuit)$ and is denoted by $\mathrm{Gal_{ant}}\left(\psi^{*},\psi_{*}\right)$ or more generally, 
\begin{equation*}
\mathrm{Gal_{ant}}\left(\psi^{*},\psi_{*}\right):(P,\diamondsuit)\to (Q,\tilde\diamondsuit) .  \tag{5.6}
\end{equation*}
\end{definition}

\begin{proposition}
Let $(P,\diamondsuit)$, $(Q,\tilde\diamondsuit)$ and $(R,\tilde{\tilde{\diamondsuit}})$ be any three binary posets. Then,
$$\mathrm{Gal\left(\psi^{*},\psi_{*}\right)}:(P,\diamondsuit)\to (Q,\tilde\diamondsuit) \land \mathrm{Gal\left(\phi^{*},\phi_{*}\right)}: (Q,\tilde{\diamondsuit})\to (R,\tilde{\tilde{\diamondsuit}})$$
\begin{equation*}
 \implies \mathrm{Gal\left(\phi^{*}\, {\circ}\, \psi^{*},\phi_{*}\, {\circ}\, \psi_{*}\right)}:(P,\diamondsuit)\to (R,\tilde{\tilde{\diamondsuit}}).   \tag{5.7}
\end{equation*}
\end{proposition}

\begin{proposition}
We have 
\begin{equation*}
 \mathrm{Gal\left(\psi^{*},\psi^{1}\right)}:(P,\diamondsuit)\to (Q,\tilde\diamondsuit) \land \mathrm{Gal\left(\psi^{*},\psi^{2}\right)}:(P,\diamondsuit)\to (Q,\tilde\diamondsuit)\implies \psi^{1}=\psi^{2}   \tag{5.8}
\end{equation*}
and
\begin{equation*}
\mathrm{Gal\left(\psi_{1},\psi_{*}\right)}:(P,\diamondsuit)\to (Q,\tilde\diamondsuit)\land\mathrm{Gal\left(\psi_{2},\psi_{*}\right)}:(P,\diamondsuit)\to (Q,\tilde\diamondsuit) \implies \psi_{1}=\psi_{2}.    \tag{5.9}
\end{equation*}
\end{proposition}

\begin{proposition}
Galois connections are not necessarily symmetric. Thus, 
\begin{equation*}
\mathrm{Gal\left(\psi^{*},\psi_{*}\right)}:(P,\diamondsuit)\to (Q,\tilde\diamondsuit) \nRightarrow \mathrm{Gal\left(\psi_{*},\psi^{*}\right)}:(Q,\tilde\diamondsuit)\to(P,\diamondsuit).    \tag{5.10}
\end{equation*}
\end{proposition}

Proof of the above propositions follows trivially from Def. 12 \cite{8}. \vspace{5mm}

\textbf{\textsc{Conclusion.}} In this paper, we have derived some results that would help us to study partially ordered sets under two partially ordered relations simultaneously. We call this the binary poset. Furthermore, it has been shown that binary posets exhibit the same properties as posets even under come complicated properties such as the uniqueness of maximal and minimal elements, existence of isotone, and isotone inverse. Some other properties such as isomorphism, duality principle, and Galois connections for binary posets are also derived. The derived results are applied to obtain some results on the power set. The next step in research in this direction would be to construct the concept of binary chains and most importantly, binary lattices \cite{3},\cite{4}. The Hasse diagrams will be different too. It has been known widely in the literature that Galois connections can be used to study posets under the category theory [5-8] and many other structures [9-12]. The structures that we can derive in category theory from the Galois connections of binary posets can be mathematically interesting to study. Work in this direction is in progress.

\end{document}